\documentclass{amsart}
\usepackage{amssymb,amsmath,amscd, amsfonts,latexsym}
\usepackage[top = 4cm, bottom = 4cm, left = 4cm, right = 4cm]{geometry}
\usepackage{hyperref}

\newtheorem{theorem}{Theorem}

\newtheorem{lemma}{Lemma}

\theoremstyle{remark}

\title[\title{On the  problem of Pillai with Fibonacci numbers and powers of $3$}]{On the  problem of Pillai with Fibonacci numbers and powers of $3$}

\author[Mahadi Ddamulira ]{Mahadi Ddamulira}



\address{Mahadi Ddamulira \newline
         \indent Institute of Analysis and Number Theory, Graz University of Technology \newline
        \indent Kopernikusgasse 24/II \newline
         \indent A-8010 Graz, Austria}

\email{mddamulira\char'100tugraz.at}

\begin{document}


\begin{abstract}
Consider the sequence $ \{F_{n}\}_{n\geq 0} $ of Fibonacci numbers  defined by $ F_0=0 $, $ F_1 =1$ and $ F_{n+2}=F_{n+1}+ F_{n} $ for all $ n\geq 0 $. In this paper, we find all integers $ c $ having at least two representations as a difference between a Fibonacci number and a power of $ 3 $.
\end{abstract}

\maketitle

\section{Introduction}
\noindent
We consider the sequence $ \{F_{n}\}_{n\geq 0} $ of Fibonacci numbers defined by
$$F_0=0, F_1= 1, \text{  and  } F_{n+2}= F_{n+1}+F_n \text{   for all  } n\geq 0.$$
The first few terms of the Fibonacci sequence are
$$\{F_{n}\}_{n\ge 0} = 0, 1, 1, 2, 3, 5, 8, 13, 21, 34, 55, 89, 144, 233, \ldots. $$

\noindent
In this paper, we are interested in studying the Diophantine equation
\begin{eqnarray}\label{Main1a}
F_n - 3^m = c
\end{eqnarray}
for a fixed integer $ c $ and variable integers $ n $ and $ m $. In particular, we are interested in finding those integers $ c $ admitting at least two representations as a difference between a Fibonacci number and a power of $ 3 $. This equation is a variant of the Pillai equation
\begin{eqnarray}\label{Pillai}
a^{x}-b^{y}=c
\end{eqnarray}
where $ x, y $ are non-gative integers and $ a,~b,~c $ are fixed positive integers.

\noindent
In 1936 and again in 1937, Pillai (see \cite{Pillai:1936, Pillai:1937}) conjectured that for any given integer $ c\geq 1 $, the number of positive integer solutions $ (a,b,x,y) $, with $ x\geq 2 $ and $ y\geq 2 $ to the equation \eqref{Pillai} is finite. This conjecture is still open for all $ c\neq 1 $. The case $ c=1 $ is Catalan's conjecture which was proved by Mih\u ailescu (see \cite{Mihailescu}). Pillai's work was an extension of the work of  Herschfeld (see  \cite{Herschfeld:1935, Herschfeld:1936}), who had already studied a particular case of the problem with $ (a,b)=(2,3) $. Since then, different variants of the Pillai equation have been studied. Some recent results for the different variants of the Pillai problem involving Fibonacci numbers, Tribonacci numbers, Pell numbers and the $ k $-generalized Fibonacci numbers with powers of $ 2 $ have been intesively studied in\cite{Luca16, Chim1, Chim2, Ddamulira1, Ddamulira2, Hernane, Hernandez}.

\section{Main Result}
\noindent
The main aim of this paper is to prove the following result.
\begin{theorem}\label{Main}
The only integers $ c $ having at least two representations of the form $ F_{n}-3^{m} $ are $ c\in\{ -26, -6, -1, 0, 2, 4, 7, 12\} $. Furthermore, all the representations of the above integers as $ F_{n}-3^{m} $ with integers $ n\geq 2 $ and $ m\geq 0 $ are given by
\begin{eqnarray}\label{lista1a}
-26 &=& F_{10} - 3^4 = F_2-3^3;\nonumber\\
-6 &=& F_8-3^3 = F_4 - 3^2;\nonumber\\
-1 &=& F_6-3^2 = F_3-3^1\nonumber\\
0&=& F_4-3^1 = F_2-3^0;\\
2&=&F_5-3^1 = F_4-3^0;\nonumber\\
4&=&F_7-3^2 = F_5-3^0;\nonumber\\
7&=&F_9-3^3=F_6-3^0;\nonumber\\
12&=&F_8-3^2=F_7-3^0.\nonumber
\end{eqnarray}
\end{theorem}

\section{Auxiliary results}
\noindent
In order to prove our main result Theorem \ref{Main}, we need to use several times a Baker--type lower bound for a nonzero linear form in logarithms of algebraic numbers. There are many such 
in the literature  like that of Baker and W{\"u}stholz from \cite{bawu07}.  We use the one of Matveev from \cite{MatveevII}.  
Matveev \cite{MatveevII} proved the following theorem, which is one of our main tools in this paper.

\noindent
Let $ \gamma $ be an algebraic number of degree $ d $ with minimal primitive polynomial over the integers
$$ a_{0}x^{d}+ a_{1}x^{d-1}+\cdots+a_{d} = a_{0}\prod_{i=1}^{d}(x-\gamma^{(i)}),$$
where the leading coefficient $ a_{0} $ is positive and the $ \eta^{(i)} $'s are the conjugates of $ \gamma $. Then the \textit{logarithmic height} of $ \gamma $ is given by
$$ h(\gamma) := \dfrac{1}{d}\left( \log a_{0} + \sum_{i=1}^{d}\log\left(\max\{|\gamma^{(i)}|, 1\}\right)\right).$$

In particular, if $ \gamma = p/q $ is a rational number with $ \gcd (p,q) = 1 $ and $ q>0 $, then $ h(\gamma) = \log\max\{|p|, q\} $. The following are some of the properties of the logarithmic height function $ h(\cdot) $, which will be used in the next sections of this paper without reference:
\begin{eqnarray}
h(\eta\pm \gamma) &\leq& h(\eta) +h(\gamma) +\log 2,\nonumber\\
h(\eta\gamma^{\pm 1})&\leq & h(\eta) + h(\gamma),\\
h(\eta^{s}) &=& |s|h(\eta) ~~~~~~ (s\in\mathbb{Z}). \nonumber
\end{eqnarray}

\begin{theorem}[Matveev]\label{Matveev11} Let $\gamma_1,\ldots,\gamma_t$ be positive real algebraic numbers in a real algebraic number field 
$\mathbb{K}$ of degree $D$, $b_1,\ldots,b_t$ be nonzero integers, and assume that
\begin{equation}
\label{eq:Lambda}
\Lambda:=\gamma_1^{b_1}\cdots\gamma_t^{b_t} - 1,
\end{equation}
is nonzero. Then
$$
\log |\Lambda| > -1.4\times 30^{t+3}\times t^{4.5}\times D^{2}(1+\log D)(1+\log B)A_1\cdots A_t,
$$
where
$$
B\geq\max\{|b_1|, \ldots, |b_t|\},
$$
and
$$A
_i \geq \max\{Dh(\gamma_i), |\log\gamma_i|, 0.16\},\qquad {\text{for all}}\qquad i=1,\ldots,t.
$$
\end{theorem} 
\noindent
During the course of our calculations, we get some upper bounds on our variables which are too large, thus we need to reduce them. To do so, we use some results from the theory of continued fractions. Specifically, for a nonhomogeneous linear form in two integer variables, we use a slight variation of a result due to Dujella and Peth{\H o} \cite{dujella98}, which  itself is a generalization of a result of Baker and Davenport \cite{BD69}.

\noindent
For a real number $X$, we write  $||X||:= \min\{|X-n|: n\in\mathbb{Z}\}$ for the distance from $X$ to the nearest integer.
\begin{lemma}[Dujella, Peth\H o]\label{Dujjella}
Let $M$ be a positive integer, $p/q$ be a convergent of the continued fraction of the irrational number $\tau$ such that $q>6M$, and  $A,B,\mu$ be some real numbers with $A>0$ and $B>1$. Let further 
$\varepsilon: = ||\mu q||-M||\tau q||$. If $ \varepsilon > 0 $, then there is no solution to the inequality
$$
0<|u\tau-v+\mu|<AB^{-w},
$$
in positive integers $u,v$ and $w$ with
$$ 
u\le M \quad {\text{and}}\quad w\ge \dfrac{\log(Aq/\varepsilon)}{\log B}.
$$
\end{lemma}

The above lemma cannot be applied when $\mu=0$ (since then $\varepsilon<0$). In this case, we use the following criterion of Legendre.

\begin{lemma}[Legendre]
\label{Legendre}
Let $\tau$ be real number and $x,y$ integers such that
\begin{equation}
\label{eq:continuedfraction}
\left|\tau-\frac{x}{y}\right|<\frac{1}{2y^2}.
\end{equation}
Then $x/y=p_k/q_k$ is a convergent of $\tau$. Furthermore, 
\begin{equation}
\label{eq:continuedfraction1}
\left|\tau-\frac{x}{y}\right|\ge \frac{1}{(a_{k+1}+2)y^2}.
\end{equation}
\end{lemma}
\noindent
Finally, the following lemma is also useful. It is Lemma 7 in \cite{guzmanluca}. 

\begin{lemma}[G\'uzman, Luca]
\label{gl}
If $m\geqslant 1$, $T>(4m^2)^m$  and $T>x/(\log x)^m$, then
$$
x<2^mT(\log T)^m.
$$
\end{lemma}

\section{Proof of Theorem \ref{Main}}
\noindent
Assume that there exist positive integers $ n, m, n_1, m_1 $ such that $ (n,m)\neq (n_{1}, m_{1}) $, and $$F_n -3^m = F_{n_{1}}-3^{m_1}.$$
Without loss of generality, we can assume that $ m\geq m_1 $. If $ m=m_1 $, then $ F_n =F_{n_1} $, so $ (n,m)= (n_{1}, m_{1}) $, which gives a contradiction to our assumption. Thus $ m>m_1 $. Since 
\begin{eqnarray}
F_n - F_{n_1} = 3^m -3^{m_1},\label{Fibo11a}
\end{eqnarray} 
and the right-hand side is positive, we get that the left-hand side is also positive and so $ n>n_{1} $. 

\noindent
Using the Binet formula
\begin{eqnarray}\label{Binet1a}
F_k = \dfrac{\alpha^{k}- \beta^k}{\sqrt{5}} \text{   for all  } k\geq 0,
\end{eqnarray}
\noindent
where $ (\alpha,\beta):=\left( \frac{1+\sqrt{5}}{2}, \frac{1-\sqrt{5}}{2}\right) $ are the roots of the equation $ x^2-x-1 =0 $, which is the characteristic equation of the Fibonacci sequence. One can easily prove by induction that
\begin{eqnarray}
\alpha^{k-2}\leq F_{n}\leq \alpha^{k-1} \text{     for all    } k\geq 1.
\end{eqnarray}
\noindent
Using the equation \eqref{Fibo11a}, we get
\begin{eqnarray}
&\alpha^{n-4}\leq F_{n-2}\leq F_n-F_{n_1}=3^m-3^{m_1} <3^m,\label{kkaa1}\\
&\alpha^{n-1}\geq F_{n}\geq F_n-F_{n_1}=3^m-3^{m_1} \geq 3^{m-1},
\end{eqnarray}
from which we get that
\begin{eqnarray}\label{Fiboo1b}
1+\left(\dfrac{\log 3}{\log\alpha}\right)(m-1) < n< \left(\dfrac{\log 3}{\log\alpha}\right)m+4.
\end{eqnarray}
\noindent
If $ n<300 $, then $ m\leq 200 $. We ran a \textit{Mathematica} program for $ 2\leq n_{1}<n\leq 300 $ and $ 0\leq m_1 < m\leq 200 $ and found only the solutions from the list \eqref{lista1a}. From now, we assume that $ n\geq 300 $. Note that the inequality \eqref{Fiboo1b} implies that $ 2m<n $. Therefore, to solve the Diophatine equation \eqref{Main1a}, it suffices to find an upper bound for $ n $.

\subsection{Bounding $ n $}
By substituting the Binet formula \eqref{Binet1a} in the Diophantine equation \eqref{Main1a}, we get
\begin{eqnarray*}
\left|\dfrac{\alpha^{n}}{\sqrt{5}}-3^{m}\right| &=& \left|\dfrac{\beta^{n}}{\sqrt{5}}+\dfrac{\alpha^{n_1}-\beta^{n_1}}{\sqrt{5}}-3^{m_1}\right| ~\leq ~ \dfrac{\alpha^{n_1}+2}{\sqrt{5}}+3^{m_1}\\
&\leq & \dfrac{2\alpha^{n_1}}{\sqrt{5}}+3^{m_1}~<~3\max\{\alpha^{n_1}, 3^{m_1}\}.
\end{eqnarray*}
Multiplying through by $ 3^{-m} $, using the relation \eqref{kkaa1} and using the fact that $ \alpha < 3 $, we get
\begin{eqnarray}
&\left|(\sqrt{5})^{-1}\alpha^{n}3^{-m} -1\right|~<~ 3\max\left\{\dfrac{\alpha^{n_1}}{3^{m}}, 3^{m_{1}-m}\right\}~<~\max\{\alpha^{n_1-n+5}, 3^{m_1-m+1}\}.\label{kkaa2}
\end{eqnarray}
For the left-hand side, we apply the result of Matveev, Theorem \ref{Matveev11} with the following data
\begin{eqnarray*}
t=3, ~~ \gamma_1=\sqrt{5}, ~~\gamma_2 = \alpha, ~~\gamma_3 = 3, ~~b_1=-1, ~~b_2 = n, ~~b_3=-m.
\end{eqnarray*}
Through out we work with the field $ \mathbb{K}:=\mathbb{Q}(\sqrt{5}) $ with $ D=2 $. Since $ \max\{1,n,m\}\leq 2n $, we take $ B:=2n $. Furthermore, we take $ A_1:=2h(\gamma_1) = \log 5 $, $ A_2:=2h(\gamma_2) = \log \alpha $, $ A_3:=2h(\gamma_1) = 2\log 3 $. We put
\begin{eqnarray*}
\Lambda = (\sqrt{5})^{-1}\alpha^{n}3^{-m} -1.
\end{eqnarray*} 
First we check that $ \Lambda \neq 0 $, if it were, then $ \alpha^{2n}\in\mathbb{Q} $, a contradiction. Thus, $ \Lambda \neq 0 $. Then by Matveev's theorem, the left-hand side of \eqref{kkaa2} is bounded as
\begin{eqnarray*}
\log|\Lambda|> -1.4\cdot 30^6\cdot 3^{4.5}\cdot 2^2(1+\log 2)(1+\log 2n)(\log 5)(\log\alpha)(2\log 3).
\end{eqnarray*}
By comparing with \eqref{kkaa2}, we get
\begin{eqnarray*}
\min\{(n-n_1-5)\log\alpha, (m-m_1-1)\log 3\}<1.04\times 10^{12}(1+\log 2n),
\end{eqnarray*}
which gives
\begin{eqnarray*}
\min\{(n-n_1)\log\alpha, (m-m_1)\log 3\}<1.24\times 10^{12}(1+\log 2n).
\end{eqnarray*}
Now we split the argument into two cases

\textbf{Case 1.}
$\min\{(n-n_1)\log\alpha, (m-m_1)\log 3\} = (n-n_1)\log\alpha$.

\noindent
In this case, we rewrite \eqref{Fibo11a} as
\begin{eqnarray*}
\left|\left(\dfrac{\alpha^{n}-\alpha^{n_1}}{\sqrt{5}}\right)-3^{m}\right|=\left|\left(\dfrac{\beta^{n}-\beta^{n_1}}{\sqrt{5}}\right)-3^{m_1}\right|~<~ 1+3^{m_1}~\leq ~3^{m_1+1},
\end{eqnarray*}
which implies
\begin{eqnarray}\label{kkaa3}
\left|\left(\dfrac{\alpha^{n-n_1}-1}{\sqrt{5}}\right)\alpha^{n_1}3^{-m}-1\right|&<&3^{m_1-m+1}.
\end{eqnarray}
We put
\begin{eqnarray*}
\Lambda_1 = \left(\dfrac{\alpha^{n-n_1}-1}{\sqrt{5}}\right)\alpha^{n_1}3^{-m}-1.
\end{eqnarray*}
To see that $ \Lambda_1 \neq 0$, for if $ \Lambda_1 =0$, then
\begin{eqnarray*}
\alpha^{n}-\alpha^{n_1}= \sqrt{5}\cdot 3^{m}.
\end{eqnarray*}
By conjugating the above relation in $ \mathbb{K} $, we get that
\begin{eqnarray*}
\beta^{n}-\beta^{n_1}= -\sqrt{5}\cdot 3^{m}.
\end{eqnarray*}
The absolute value of the left-hand side is at most $ |\beta^{n}-\beta^{n_1}|\leq |\beta|^{n}+|\beta|^{n_1}<2 $, while the absolute value of the right-hand side is at least $|-\sqrt{5}\cdot 3^{m}|\geq \sqrt{5}> 2$ for all $ m\geq 0 $, which is a contradiction. 

\noindent
We apply Theorem \ref{Matveev11} on the left-hand side of \eqref{kkaa3} with the data
\begin{eqnarray*}
t=3, ~~\gamma_1=\dfrac{\alpha^{n-n_1}-1}{\sqrt{5}},~~ \gamma_2=\alpha, ~~\gamma_3=3, ~~ b_1=1, ~~b_2=n_1, ~~b_3=-m.
\end{eqnarray*}
The minimal polynomial of $ \gamma_1 $ divides
\begin{eqnarray*}
5X^2-5F_{n-n_1}X-((-1)^{n-n_1}+1-L_{n-n_1}),
\end{eqnarray*}
where $ \{L_{k}\}_{k\geq 0} $ is the Lucas companion sequence of the Fibonacci sequence given by $ L_0=2, ~~L_1=1,~~L_{k+2}=2L_{k+1}+L_k $ for all $ k\geq 0 $, for which the Binet formula for its general term is given by
\begin{eqnarray*}
L_k=\alpha^{k}+\beta^{k} ~~\text{  for all  } k\geq 0.
\end{eqnarray*}
Thus, we obtain
\begin{eqnarray}
h(\gamma_1)&\leq& \dfrac{1}{2}\left(\log 5+\log\left(\dfrac{\alpha^{n-n_1}+1}{\sqrt{5}}\right)\right)~<~\dfrac{1}{2}\log(4\sqrt{5}\alpha^{n-n_1})\nonumber\\
&<&\dfrac{1}{2}(n-n_1+4)\log\alpha~<~8.4\times 10^{11}(1+\log 2n).\label{height1a}
\end{eqnarray}
So, we can take $ A_1:=16.8\times 10^{11}(1+\log 2n) $. Furthermore, as before, we take $ A_2:=\log\alpha $ and $ A_3:=2\log 3 $. Finally, since $ \max\{1, n_1, m\}\leq 2n $, we can take $ B:=2n $. Then, we get
\begin{eqnarray*}
\log|\Lambda_1|> -1.4\cdot 30^6\cdot 3^{4.5}\cdot 2^2(1+\log 2)(1+\log 2n)(16.8\times 10^{11}(1+\log 2n))(\log\alpha)(2\log 3).
\end{eqnarray*}
Then,
\begin{eqnarray*}
\log|\Lambda_1|>-1.72\times 10^{24}(1+\log 2n)^2.
\end{eqnarray*}
By comparing the above relation with \eqref{kkaa3}, we get that
\begin{eqnarray}
(m-m_1)\log 3 < 1.80\times 10^{24}(1+\log 2n)^2.
\end{eqnarray}

\textbf{Case 2.}
$\min\{(n-n_1)\log\alpha, (m-m_1)\log 3\} =  (m-m_1)\log 3$.

\noindent
In this case, we rewrite \eqref{Fibo11a} as
\begin{eqnarray*}
\left|\dfrac{\alpha^{n}}{\sqrt{5}}-(3^{m-m_1}-1)\cdot 3^{m_1}\right|=\left|\dfrac{\beta^n+\alpha^{n_1}-\beta^{n_1}}{\sqrt{5}}\right|~<~\dfrac{\alpha^{n_1}+2}{\sqrt{5}}~<~\alpha^{n_1},
\end{eqnarray*}
which implies that
\begin{eqnarray}\label{kkaa4}
|(\sqrt{5}(3^{m-m_1}-1))^{-1}\alpha^{n}3^{-m_1}-1|&<&\dfrac{\alpha^{n_1}}{3^{m}-3^{m_1}}~\leq ~\dfrac{3\alpha^{n_1}}{3^{m}}\nonumber\\&<&3\alpha^{n_1-n+4}~<~\alpha^{n_1-n+5}.
\end{eqnarray}
We put
\begin{eqnarray*}
\Lambda_{2}=(\sqrt{5}(3^{m-m_1}-1))^{-1}\alpha^{n}3^{-m_1}-1.
\end{eqnarray*}
Clearly, $ \Lambda_2\neq 0 $, for if $ \Lambda_2 =0 $, then $ \alpha^{2n}\in\mathbb{Q} $, which is a contradiction. We again apply Theorem \ref{Matveev11} with the following data
\begin{eqnarray*}
t=3,~~ \gamma_{1}=\sqrt{5}(3^{m-m_1}-1), ~~\gamma_2=\alpha, ~~\gamma_3=\alpha, ~~b_1=-1, ~~b_2=n, ~~b_3= -m_1.
\end{eqnarray*}
The minimal polynomial of $ \gamma_1 $ is $ X^2-5(3^{m-m_1}-1)^2 $. Thus,
\begin{eqnarray*}
h(\gamma_{1})=\log\left(\sqrt{5}(3^{m-m_1}-1)\right)~<~(m-m_1+1)\log 3~ <~ 1.25\times 10^{12}(1+\log 2n).
\end{eqnarray*}
So, we can take $ A_1:=2.5\times 10^{12}(1+\log 2n) $. Further, as in the previous applications, we take $ A_2:=\log\alpha $ and $ A_3:=2\log 3 $. Finally, since $ \max\{1, n, m_1\}\leq 2n $, we can take $ B:=2n $. Then, we get
\begin{eqnarray*}
\log|\Lambda_2|> -1.4\cdot 30^6\cdot 3^{4.5}\cdot 2^2(1+\log2)(1+\log 2n)(2.5\times 10^{12}(1+\log2n))(\log\alpha)(2\log 3).
\end{eqnarray*}
Thus,
\begin{eqnarray*}
\log|A_2|>- 2.56\times 10^{24}(1+\log 2n)^2.
\end{eqnarray*}
Now, by comparing with \eqref{kkaa4}, we get that 
\begin{eqnarray}
(n-n_1)\log\alpha < 2.58\times 10^{24}(1+\log 2n)^2.
\end{eqnarray}
Therefore, in both Case 1 and Case 2, we have
\begin{eqnarray}
\min\{(n-n_1)\log\alpha, (m-m_1)\log 3\}&<&1.24\times 10^{12}(1+\log 2n),\nonumber\\
\max\{(n-n_1)\log\alpha, (m-m_1)\log 3\}&<&2.58\times 10^{24}(1+\log 2n)^2.\label{height1b}
\end{eqnarray}

\noindent
Finally, we rewrite the equation \eqref{Fibo11a} as
\begin{eqnarray*}
\left|\dfrac{(\alpha^{n-n_1}-1)}{\sqrt{5}}\alpha^{n_1}-(3^{m-m_1}-1)\cdot 3^{m_1}\right|=\left|\dfrac{\beta^{n}-\beta^{n_1}}{\sqrt{5}}\right|~<~|\beta|^{n_1}=\alpha^{-n_1}.
\end{eqnarray*}
Dividing through by $ 3^m-3^{m_1} $, we get
\begin{eqnarray}
\left|\left(\dfrac{\alpha^{n-n_1}-1}{\sqrt{5}(3^{m-m_1}-1)}\right)\alpha^{n_1}3^{-m_1}-1\right|&<&\dfrac{1}{\alpha^{n_1}(3^{m}-3^{m_1})}\leq \dfrac{3}{\alpha^{n_1}\cdot 3^{m}}\nonumber\\
&\leq& 3\alpha^{-(n+n_1-4)}~\leq ~\alpha^{4-n},\label{kkaa5}
\end{eqnarray}
since $ 3< \alpha \leq \alpha^{n_1} $.  We again apply Theorem \ref{Matveev11} on the left-hand side of \eqref{kkaa5} with the data
\begin{eqnarray*}
t=3, ~~\gamma_{1}=\dfrac{\alpha^{n-n_1}-1}{\sqrt{5}(3^{m-m_1}-1)}, ~~\gamma_2=\alpha, ~~\gamma_3=3, ~~b_1=1, ~~b_2=n_1, ~~b_3=-m_1.
\end{eqnarray*}
By using the algebraic properties of the logarithmic height function, we get
\begin{eqnarray*}
h(\gamma_1) &=& h\left(\dfrac{\alpha^{n-n_1}-1}{\sqrt{5}(3^{m-m_1}-1)}\right) ~\leq~h\left(\dfrac{\alpha^{n-n_1}-1}{\sqrt{5}}\right)+h(3^{m-m_1}-1)\\
&<&\dfrac{1}{2}(n-n_1+4)\log\alpha+(m-m_1)\log 3~<2.80\times 10^{24}(1+\log 2n)^2,
\end{eqnarray*}
where in the above inequalities, we used the argument from \eqref{height1a} as well as the bounds \eqref{height1b}. Thus, we can take $ A_1:=5.60\times 10^{24}(1+\log 2n) $, and again as before $ A_2:=\log\alpha $ and $ A_3:=2\log 3 $. If we put
\begin{eqnarray*}
\Lambda_3 = \left(\dfrac{\alpha^{n-n_1}-1}{\sqrt{5}(3^{m-m_1}-1)}\right)\alpha^{n_1}3^{-m_1}-1,
\end{eqnarray*}
we need to show that $ \Lambda_3\neq 0 $. If not, $ \Lambda_3=0 $ leads to
\begin{eqnarray*}
\alpha^{n}-\alpha^{n_1}=\sqrt{5}(3^{m}-3^{m_1}).
\end{eqnarray*}
A contradiction is reached upon a conjuagtion in $ \mathbb{K} $ and by taking absolute values on both sides. Thus, $ \Lambda_3 \neq 0$. Applying Theorem \ref{Matveev11} gives
\begin{eqnarray*}
\log|\Lambda_3|>-1.4\cdot 30^6 \cdot 3^{4.5}\cdot 2^2 (1+\log2)(1+\log 2n)(5.6\times 10^{24}(1+\log 2n)^2)(\log\alpha)(2\log 3),
\end{eqnarray*}
a comparison with \eqref{kkaa5} gives
\begin{eqnarray*}
(n-4)&<&3\times 10^{36}(1+\log 2n)^3,
\end{eqnarray*}
or 
\begin{eqnarray}\label{fala1a}
2n&<&6.2\times 10^{36}(1+\log 2n)^3.
\end{eqnarray}
Now by applying Lemma \ref{gl} on \eqref{fala1a} with the data $ m=3 $, $ T=6.2\times 10^{36} $ and $ x=2n $, leads to $ n< 2\times 10^{40} $.

\subsection{Reducing the bound for $ n $}
\noindent
We need to reduce the above bound for $ n $ and to do so we make use of Lemma \ref{Dujjella} several times. To begin, we return to \eqref{kkaa2} and put
\begin{eqnarray*}
\Gamma:=n\log\alpha -m\log 3 -\log(\sqrt{5}).
\end{eqnarray*}
For technical reasons we assume that $ \min\{n-n_1, m-m_1\}\geq 20 $. We go back to the inequalities for $ \Lambda $, $ \Lambda_1 $ and  $ \Lambda_2 $, Since we assume that  $ \min\{n-n_1, m-m_1\}\geq 20 $ we get $ |e^{\Gamma}-1|=|\Lambda|< \frac{1}{4} $. Hence, $ |\Lambda|<\frac{1}{2} $ and since the inequality $ |y|<2|e^{y}-1| $ holds for all $ y\in\left(-\frac{1}{2}, \frac{1}{2}\right) $, we get
\begin{eqnarray*}
|\Gamma|<2\max\{\alpha^{n_1-n+5}, 3^{m_1-m+1}\}\leq \max\{\alpha^{n_1-n+6}, 3^{m_1-m+2}\}.
\end{eqnarray*}
Assume that $ \Gamma >0 $. We then have the inequality
\begin{eqnarray*}
0<n\left(\dfrac{\log\alpha}{\log3}\right)-m+\dfrac{\log(1/\sqrt{5})}{\log 3}&<&\max\left\{\dfrac{\alpha^{8}}{(\log 3)\alpha^{n-n_1}}, \dfrac{6}{(\log3)3^{m-m_1}}\right\}.\\
&<&\max\{45\alpha^{-(n-n_1)}, 8\cdot 3^{-(m-m_1)}\}.
\end{eqnarray*}
We apply Lemma \ref{Dujjella} with the data
\begin{eqnarray*}
\tau = \dfrac{\log\alpha}{\log 3}, ~~~\mu = \dfrac{\log(1/\sqrt{5})}{\log 3}, ~~~(A,B)=(45, \alpha) ~~\text{  or } (8,3).
\end{eqnarray*}
Let $ \tau = [a_{0}; a_{1}, a_{2}, \ldots]=[0; 2, 3, 1, 1, 6, 1, 49, 1, 2, 2, 1, 1, 2, 1, 2, 2, 1, 10, 3, \ldots] $ be the continued fraction of $ \tau $. We choose $M:=2\times 10^{40}$ and consider the 91-th convergent
\begin{eqnarray*}
\dfrac{p}{q}=\dfrac{p_{91}}{q_{91}}=\dfrac{487624200385184167130255744232737921512174859336581}{1113251817385764505972408650620147577750763395186265}.
\end{eqnarray*}
It satisfies $ q=q_{91}>6M $. Furthermore, it yields $ \varepsilon > 0.50$, and therefore either
\begin{eqnarray*}
n-n_1\leq \dfrac{\log(45q/\varepsilon)}{\log\alpha} < 254, ~~\text{  or  } ~~ m-m_1\leq \dfrac{\log (8q/\varepsilon)}{\log 3}<110.
\end{eqnarray*}

\noindent
In the case $ \Gamma <0 $, we consider the inequality
\begin{eqnarray*}
m\left(\dfrac{\log 3}{\log\alpha}\right)-n+\dfrac{\log (\sqrt{5})}{\log\alpha}&<&\max\left\{\dfrac{\alpha^{8}}{\log\alpha}\alpha^{-(n-n_1)}, ~\dfrac{8}{\log\alpha}\cdot 3^{-(m-m_1)}\right\}\\
&<&\max\{98\alpha^{-(n-n_1)}, ~18\cdot 3^{-(m-m_1)}\}.
\end{eqnarray*}
We then apply Lemma \ref{Dujjella} with the data
\begin{eqnarray*}
\tau = \dfrac{\log 3}{\log\alpha}, ~~ \mu = \dfrac{\log \sqrt{5}}{\log\alpha}, ~~ (A,B)=(98, \alpha), ~~\text{  or  } ~ (18,3).
\end{eqnarray*}
Let $ \tau = [a_0; a_1, a_2,  \ldots] = [2; 3, 1, 1, 6, 1, 49, 1, 2, 2, 1, 1, 2, 1, 2, 2, 1, 10, 3, 12, \ldots] $ be the continued fraction of $ \tau $. Again, we choose $ M=2\times 10^{40} $, and in this case we consider the 101-th convergent
\begin{eqnarray*}
\dfrac{p}{q}=\dfrac{p_{101}}{q_{101}}=\dfrac{106360048375891410642967692492903700137161881169662}{56228858848524361385900581302251812795713192394033},
\end{eqnarray*}
which satisfies $ q=q_{101}>6M $. Further, this yields $ \varepsilon>0.125 $, and therefore either
\begin{eqnarray*}
n-n_1\leq \dfrac{\log(98q/\varepsilon)}{\log\alpha}<254 ~~~, \text{  or  }~~ m-m_1\leq \dfrac{\log(18q/\varepsilon)}{\log 3}< 110.
\end{eqnarray*}
These bounds agree with the bounds obtained in the case $ \Gamma>0 $. As a conclusion, we have that either $ n-n_1\leq 253 $ or $ m-m_1 \leq 109$ whenever $ \Gamma\neq 0 $.

\noindent
Now, we distinguish between the cases $ n-n_1\leq 253 $ and $ m-m_1\leq 109 $. First, we assume that $ n-n_1\leq 253 $. In this case we consider the inequality for $ \Lambda_1 $, \eqref{kkaa3} and also assume that $ m-m_1\leq 20 $. We put
\begin{eqnarray*}
\Gamma_1=n_1\log\alpha-m\log3+\log\left(\dfrac{\alpha^{n-n_1}}{\sqrt{5}}\right).
\end{eqnarray*} 
Then inequality \eqref{kkaa3} implies that
\begin{eqnarray*}
|\Gamma_1|<\dfrac{6}{3^{m-m_1}}.
\end{eqnarray*}
If we further assume that $ \Gamma_1>0 $, we then get
\begin{eqnarray*}
0<n_{1}\left(\dfrac{\log\alpha}{\log 3}\right)-m+\dfrac{\log((\alpha^{n-n_1}-1)/\sqrt{5})}{\log3}<\dfrac{6}{(\log 3)3^{m-m_1}}<\frac{6}{3^{m-m_1}}.
\end{eqnarray*}
Again we apply Lemma \ref{Dujjella} with the same $ \tau $ as in the case $ \Gamma>0 $. We use the 91-th convergent $ p/q=p_{91}/q_{91} $ of $ \tau $ as before. But in this case we choose $ (A,B):=(8,3) $ and use
\begin{eqnarray*}
\mu_l=\dfrac{\log((\alpha^{l}-1)/\sqrt{5})}{\log 3},
\end{eqnarray*}
instead of $ \mu $ for each possible value of $ l:=n-n_1 \in [1,2,\ldots, 253] $. We have problems at $ l \in\{4,12\} $. We discard these values for now and we will treat them later. For the remaining values of $ l $, we get $ \varepsilon>0.0005 $. Hence by Lemma \ref{Dujjella}, we get
\begin{eqnarray*}
m-m_1~<~\dfrac{\log(8q/0.0005)}{\log 3}~<~116.
\end{eqnarray*}
Thus, $ n-n_1\leq 253 $ implies that $ m-m_1\leq 115 $, unless $ n-n_1\in\{4,12\} $. A similar conclusion is reached when $ \Gamma_1<0 $ with the same two exceptions for $ n-n_1\in\{4,12\} $. The reason we have a problem at $ l\in\{4, 12\} $ is because
\begin{eqnarray*}
\dfrac{\alpha^{4}-1}{\sqrt{5}}=\alpha^{2}, ~~~\text{  and  }~~\dfrac{\alpha^{12}-1}{\sqrt{5}}=2^{3}\alpha^{6}.
\end{eqnarray*}
So, $ \Gamma_1=(n_1+2)\log\alpha-m\log 3 $ , or $ (n_1+6)\log\alpha-(m-3)\log3 $ when $ l=4,12 $, respectively. Thus we get that
\begin{eqnarray*}
\left|\tau - \dfrac{m}{n_1+2}\right|<\dfrac{6}{3^{m-m_1}(n_{1}+2)}, ~~\text{  or  } ~~ \left|\tau - \dfrac{m-3}{n_1+6}\right|<\dfrac{6}{3^{m-m_1}(n_{1}+6)},
\end{eqnarray*}
respectively. We assume that $ m-m_1>150 $. Then $ 3^{m-m_1}>8\times(4\times 10^{40})>8\times(n_1+6) $, therefore
\begin{eqnarray*}
\dfrac{6}{3^{m-m_1}(n_1+2)}<\dfrac{1}{3(n_1+2)^{2}}, ~~~\text{  and  } ~~ \dfrac{6}{3^{m-m_1}(n_1+6)}<\dfrac{1}{3(n_1+6)^{2}}.
\end{eqnarray*}
By Lemma \ref{Legendre}, it follows that $ m/(n_1+2) $ or $ (m-3)/(n_1+6) $ are convergents of $ \tau $, respectively. So, say one of $ m/(n_1+2) $ or $ (m-3)/(n_1+6) $ is of the form $ p_k/q_k $ for some $ k=0, 1, 2, \ldots, 92 $. Here, we use that $ q_{92} >4\times 10^{40}> n+1+6$. Then
\begin{eqnarray*}
\dfrac{1}{(a_{k}+2)q_k^2}<\left|\tau -\dfrac{p_k}{q_k}\right|.
\end{eqnarray*}
Since $ \max\{a_{k}:k=0,1,2, \ldots,92\}= 140 $, we get
\begin{eqnarray*}
\dfrac{1}{142q_k^2}<\dfrac{6}{3^{m-m_1}q_k}~~~\text{ and }~~ q_{k}~~ \text{  divides one of }~~\{n_1+2, n_1+6\}.
\end{eqnarray*}
Thus, we get
\begin{eqnarray*}
3^{m-m_1}\leq 6\times 142(n_1+6)< 6\times 142\times 4\times 10^{40},
\end{eqnarray*}
giving $ m-m_1\leq 92 $.

\noindent
Now let us turn to the case $ m-m_1\leq 109 $ and we consider the inequlity for $ \Lambda_2 $, \eqref{kkaa4}. We put
\begin{eqnarray*}
\Gamma_{2}=n\log\alpha-m_1\log 3+\log(1/(\sqrt{5}(3^{m-m_1}-1))),
\end{eqnarray*}
and we also assume that $ n-n_1\geq 20 $. We then have
\begin{eqnarray*}
|\Gamma_{2}|<\dfrac{2\alpha^{8}}{\alpha^{n-n_1}}.
\end{eqnarray*}
We assume that $ \Gamma_{2} $, then we get
\begin{eqnarray*}
0~<~n\left(\dfrac{\log\alpha}{\log 3}\right)-m_1+\dfrac{\log(1/(\sqrt{5}(3^{m-m_1}-1))}{\log\alpha}~<~\dfrac{3\alpha^{8}}{(\log 3)\alpha^{n-n_1}}~<~\dfrac{130}{\alpha^{n-n_1}}.
\end{eqnarray*}
We apply again Lemma \ref{Dujjella} with the same $ \tau, ~ q, ~ M, ~~ (A,B):=(130, \alpha) $ and 
\begin{eqnarray*}
\mu_l=\dfrac{\log(1/(\sqrt{5}(3^{l}-1)))}{\log 3}~~\text{   for   } k=1, 2, \ldots, 109.
\end{eqnarray*}
We get $ \varepsilon> 0.004 $, therefore
\begin{eqnarray*}
n-n_1<\dfrac{\log(130q/\varepsilon)}{\log\alpha}<266.
\end{eqnarray*}
A similar conclusion is reached when $ \Gamma_2<0 $. To conclude, we first get that either $ n-n_1\leq 253 $ or $ m-m_1\leq 109 $. If $ n-n_1\leq 253 $, then $ m-m_1\leq 115 $, and if $ m-m_1\leq 109 $ then $ n-n_1\leq 265 $. Thus, we conclude that we always have $ n-n_1\leq 265 $ and $ m-m_1\leq 115 $.

\noindent
Finally we go to the inequality of $ \Lambda_{3} $, \eqref{kkaa5}. We put
\begin{eqnarray*}
\Gamma_3=n_1\log\alpha-m_1\log 3+\log\left(\dfrac{\alpha^{n-n_1}-1)}{\sqrt{5}(3^{m-m_1}-1)}\right).
\end{eqnarray*}
Since $ n\geq 300 $, the inequality \eqref{kkaa5} implies that
\begin{eqnarray*}
|\Gamma_{3}|<\dfrac{3}{\alpha^{n-4}}=\dfrac{3\alpha^{4}}{\alpha^{n}}.
\end{eqnarray*}
Assuming that $ \Gamma_3>0 $, then
\begin{eqnarray*}
0<n_1\left(\dfrac{\log\alpha}{\log 3}\right)-m_1+\dfrac{\log((\alpha^{k}-1)/(\sqrt{5}(3^{l}-1))}{\log 3}<\dfrac{3\alpha^{4}}{(\log 3)\alpha^{n}}<\dfrac{20}{\alpha^{n}},
\end{eqnarray*}
where $ (k,l):=(n-n_1, m-m_1)  $. We again apply Lemma \ref{Dujjella} with the same $ \tau, ~~q, ~~M, ~~(A,B):=(20, \alpha) $ and 
\begin{eqnarray*}
\mu_{k,l}=\dfrac{\log((\alpha^{k}-1)/(\sqrt{5}(3^{l}-1))}{\log 3}~~\text{ for } ~~ 1\leq k\leq 265, ~~1\leq l\leq 115.
\end{eqnarray*}
As before, we have a problem at $ (k,l):=(4,1), ~(12,1), ~ (8,2) $. The cases $ (k,l):=(4,1), (12,1) $ were treated before in the case of $ \Gamma_1 $. The case $ (k,l):=(8,2) $ arises because
\begin{eqnarray*}
\dfrac{\alpha^{8}-1}{\sqrt{5}(3^{2}-1)}=\dfrac{3}{8}\alpha^{4},
\end{eqnarray*}
we therefore discard the cases $ (k,l):=(4,1), ~(12,1), ~ (8,2) $ for some time. For the remaining cases, we get $ \varepsilon> 0.0015 $, so we obtain
\begin{eqnarray*}
n\leq \dfrac{\log(20q/\varepsilon)}{\log\alpha}<264.
\end{eqnarray*}
A similar conclusion is reached when $ \Gamma_{3}<0 $. Hence, $n<300$.
Now we look at the cases $ (k,l):=(4,1), ~(12,1), ~ (8,2) $. The cases $ (k,l):=(4,1), ~(12,1) $ can be treated as before when we showed that $ n-n_1 \leq 263$ implies $ m-m_1\leq 115 $. The case when $ (k,l)=(8,2) $ can be delt with in a similar way. Namely, it gives that
\begin{eqnarray*}
|(n_1+4)\tau - m_1|<\dfrac{20}{\alpha^{n}}.
\end{eqnarray*}
Therefore,
\begin{eqnarray}\label{kkaa6}
\left|\tau - \dfrac{m_1}{n_1+4}\right|<\dfrac{20}{(n_1+4)\alpha^{n}}.
\end{eqnarray}
Since $ n\geq 300 $, we have $ \alpha^{n}>2\times 20\times (4\times 10^{40})>40(n_1+4) $. This shows that the right hand side of the above inequality, \eqref{kkaa6} is at most $ 2/(n_1+4)^2 $. By Lemma \ref{Legendre}, we get that $ m_1/(n_1+4)=p_k/q_k $ for some $ k=1,2, \ldots, 92 $. We then get by a similar argument as before that
\begin{eqnarray*}
\alpha^{n}<20\times 142\times (4\times 10^{40}),
\end{eqnarray*}
which gives $ n\leq 211 $. Therefore, the conclusion is that $ n<300 $ holds also in the case $ (k,l)=(8,2) $. However, this contradicts our working assumption that $ n\geq 300 $. This completes the proof of Theorem \ref{Main}.

\section*{acknowledgements}
\noindent
The author was supported by the Austrian Science Fund (FWF) grants: F5510-N26 -- Part of the special research program (SFB), ``Quasi Monte Carlo Metods: Theory and Applications'', P26114-N26 --``Diophantine Problems: Analytic, geometric and computational aspects'' and W1230 --``Doctoral Program Discrete Mathematics''.

 \end{document}